\def\version{}

\documentclass[reqno,twoside,11pt]{amsart}
\usepackage{cite}
\usepackage{amsmath}
\usepackage{amsfonts}
\usepackage{amssymb}
\usepackage{epsfig}
\usepackage{verbatim}
\usepackage{color}

\IfFileExists{epsf.def}{\input epsf.def}{\usepackage{epsf}}
\IfFileExists{mathptmx.sty}{\usepackage{mathptmx}}{\usepackage{mathpazo}}
\usepackage{mathrsfs}

\DeclareFontFamily{OT1}{eusb}{} \DeclareFontShape{OT1}{eusb}{m}{n} {<5> <6> <7> <8> <9> <10> <11> <12> <14.4> eusb10}{}
\DeclareMathAlphabet{\eusb}{OT1}{eusb}{m}{n}

\DeclareFontFamily{OT1}{eusm}{} \DeclareFontShape{OT1}{eusm}{m}{n} {<5> <6> <7> <8> <9> <10> <11> <12> <14.4> eusm10}{}
\DeclareMathAlphabet{\eusm}{OT1}{eusm}{m}{n}

\DeclareFontFamily{OT1}{eufm}{} \DeclareFontShape{OT1}{eufm}{m}{n} {<5> <6> <7> <8> <9> <10> <11> <12> <14.4> eufm10}{}
\DeclareMathAlphabet{\mathfrak}{OT1}{eufm}{m}{n}

\DeclareFontFamily{OT1}{fraktura}{}
\DeclareFontShape{OT1}{fraktura}{m}{n} {<5> <6> <7> <8> <9> <10> <11> <12> <13> <14.4> [1.1] eufm10}{}
\DeclareMathAlphabet{\fraktura}{OT1}{fraktura}{m}{n}

\DeclareFontFamily{OT1}{cmfi}{} \DeclareFontShape{OT1}{cmfi}{m}{n} {<5> <6> <7> <8> <9> <10> <11> <12> <13> <14.4> [0.9] cmfi10}{}
\DeclareMathAlphabet{\cmfi}{OT1}{cmfi}{b}{n}

\DeclareFontFamily{OT1}{cmss}{} \DeclareFontShape{OT1}{cmss}{m}{n} {<5> <6> <7> <8> <9> <10> <11> <12> <13> <14.4> cmss10}{}
\DeclareMathAlphabet{\cmss}{OT1}{cmss}{m}{n}

\newtheoremstyle{thm}{1.5ex}{1.5ex}{\itshape\rmfamily}{} {\bfseries\rmfamily}{}{2ex}{}

\newtheoremstyle{def}{1.5ex}{1.5ex}{\rmfamily\sl}{} {\bfseries\rmfamily}{}{2ex}{}

\newtheoremstyle{rem}{1.3ex}{1.3ex}{\rmfamily}{} {\bfseries\rmfamily}{}{2ex}{}

\newtheoremstyle{ass}{1.5ex}{1.5ex}{\rmfamily\sl}{} {\bfseries\rmfamily}{}{2ex}{}

\newenvironment{proofsect}[1] {\vskip0.1cm\noindent{\rmfamily\itshape#1.}}{\qed\vspace{0.15cm}}

\theoremstyle{thm}
\newtheorem{theorem}{Theorem}[section]
\newtheorem{lemma}[theorem]{Lemma}
\newtheorem{proposition}[theorem]{Proposition}

\newtheorem*{Main Theorem}{Main Theorem.}
\newtheorem{corollary}[theorem]{Corollary}
\newtheorem{conjecture}[theorem]{Conjecture}

\newtheorem{observation}[theorem]{Observation}

\theoremstyle{rem}
\newtheorem{remark}[theorem]{{Remark}}

\numberwithin{equation}{section}

\renewcommand{\section}{\secdef\sct\sect}
\newcommand{\sct}[2][default]{\refstepcounter{section}
\addcontentsline{toc}{section}
{{\tocsection {}{\thesection}{\!\!\!\!#1\dotfill}}{}}
\vspace{0.7cm}
\centerline{ 
\scshape\arabic{section}.\ #1} \nopagebreak \vspace{0.2cm}}
\newcommand{\sect}[1]{
\vspace{0.4cm} \centerline{\large\scshape\rmfamily #1}
\vspace{0.2cm}}

\renewcommand{\subsection}{\secdef\subsct\sbsect}
\newcommand{\subsct}[2][default]{\refstepcounter{subsection}
\addcontentsline{toc}{subsection}
{{\tocsection{\!\!}{\hspace{1.2em}\thesubsection}{\!\!\!\!#1\dotfill}}{}}
\nopagebreak\vspace{0.45\baselineskip} {\flushleft\bf
\thesection.\arabic{subsection}~\bf #1.~}
\\*[3mm]\noindent
\nopagebreak}
\newcommand{\sbsect}[1]{
\vspace{0.1cm}\noindent
\textbf{#1.~}\vspace{0.1cm}}

\renewcommand{\subsubsection}{%
\secdef \subsubsect\sbsbsect}
\newcommand{\subsubsect}[2][default]{%
\refstepcounter{subsubsection} 
\addcontentsline{toc}{subsubsection}{{\tocsection{\!\!}
{\hspace{3.05em}\thesubsubsection}{\!\!\!\!#1\dotfill}}{}}
\nopagebreak
\vspace{0.15\baselineskip} \nopagebreak {\flushleft\rmfamily
\itshape\arabic{section}.\arabic{subsection}.\arabic{subsubsection}
\ \rmfamily #1\/.}\ }
\newcommand{\sbsbsect}[1]{\vspace{0.1cm}\noindent
\rmfamily \itshape
\arabic{section}.\arabic{subsection}.\arabic{subsubsection} \
\sffamily #1\/.\ }

\renewcommand{\caption}[1]{%
\vglue0.5cm
\refstepcounter{figure}
\begin{center}
\begin{minipage}[c]{0.8\textwidth}\small {\sc Fig.~\thefigure\ }#1\end{minipage}
\end{center}
}

\newcommand{\textd}{\text{\rm d}\mkern0.5mu}

\newcommand{\texte}{\text{\rm  e}\mkern0.7mu}

\newcommand{\Var}{\text{\rm Var}}
\newcommand{\1}{{1\mkern-4.5mu\textrm{l}}}
\renewcommand{\1}{\text{\sf 1}}

\newcommand{\II}{\mathcal I}

\newcommand{\N}{\mathbb N}

\newcommand{\R}{\mathbb R}

\newcommand{\Z}{\mathbb Z}

\newcommand{\twoeqref}[2]{(\ref{#1}--\ref{#2})}

\newcommand{\cc}{{\text{\rm c}}}

\def\myffrac#1#2 in #3{\raise 2.6pt\hbox{$#3 #1$}\mkern-1.5mu\raise 0.8pt\hbox{$#3/$}\mkern-1.1mu\lower 1.5pt\hbox{$#3 #2$}}

\newcommand{\wt}{\widetilde}

\newcommand{\laweq}{\,\overset{\text{\rm law}}=\,}

\newcommand\independent{\protect\mathpalette{\protect\independenT}{\perp}}
\def\independenT#1#2{\mathrel{\rlap{$#1#2$}\mkern3mu{#1#2}}}

\newcommand{\Sym}{\text{\rm Sym}}

\begin{document}

\title[Distance in long-range percolation \hfill \version\hfill]
{\large Sharp asymptotic for the chemical distance\\in long-range percolation}

\author[\hfill  \version \hfill Biskup and Lin]
{Marek~Biskup$^{1,2}$ and Jeffrey Lin$^1$}
\thanks{\hglue-4.5mm\fontsize{9.6}{9.6}\selectfont\copyright\,\textrm{2017}\ \ \textrm{M.~Biskup, J. Lin.
Reproduction, by any means, of the entire
article for non-commercial purposes is permitted without charge.\vspace{2mm}}}
\maketitle

\vspace{-5mm}
\centerline{\textit{
$^1$Department of Mathematics, UCLA, Los Angeles, California, USA}}
\centerline{\textit{
$^2$Center for Theoretical Study, Charles University, Prague, Czech Republic}}


\vskip0.5cm
\begin{quote}
\footnotesize \textbf{Abstract:}
We consider instances of long-range percolation on~$\Z^d$ and~$\R^d$, where points at distance~$r$ get connected by an edge with probability proportional to~$r^{-s}$, for $s\in (d,2d)$, and study the asymptotic of the graph-theoretical (a.k.a.~chemical) distance $D(x,y)$ between~$x$ and~$y$ in the limit as~$|x-y|\to\infty$. For the model on~$\Z^d$  we show that, in probability as~$|x|\to\infty$, the distance $D(0,x)$ is squeezed between two positive multiples of~$(\log r)^\Delta$, where $\Delta:=1/\log_2(1/\gamma)$ for~$\gamma:=s/(2d)$. For the model on~$\R^d$  we show that $D(0,xr)$ is, in probability as~$r\to\infty$ for any nonzero~$x\in\R^d$, asymptotic to $\phi(r)(\log r)^\Delta$ for~$\phi$ a positive, continuous (deterministic) function obeying $\phi(r^\gamma)=\phi(r)$ for all~$r>1$. The proof of the asymptotic scaling is based on a subadditive argument along a continuum  of doubly-exponential sequences of scales. The results strengthen considerably the conclusions obtained earlier by the first author. Still, significant open questions remain.
\end{quote}

\newcommand{\Ddis}{D^{\text{\,\rm dis}}}

\section{Introduction}
\nopagebreak\vglue-4mm
\subsection{The model and main results}
\nopagebreak\noindent
Long-range percolation is a tool to expand connectivity of a given graph by adding, at random, edges between far-away vertices. Although arising from questions in mathematical physics (Dyson~\cite{Dyson}, Fr\"ohlich and Spencer~\cite{FS}), the problem was recognized quickly to pose interesting challenges for probability (Schulman~\cite{Schulman}, Newman and Schulman~\cite{NS}, Aizenman and Newman~\cite{AN}, Aizenman, Kesten and Newman~\cite{AKN}). More recently, instances of long-range percolation have been used as an ambient medium for other stochastic processes (e.g., Berger~\cite{Berger-RW}, Benjamini, Berger and Yadin~\cite{BBY}, Crawford and Sly~\cite{CS1,CS2}, Misumi~\cite{Misumi}, Kumagai and Misumi~\cite{KM}). The overarching theme here is the geometry of random networks.

In this paper we consider two models of long-range percolation on~$\R^d$. One of these is set on the hypercubic lattice~$\Z^d$ (endowed, \textit{a priori}, with its nearest-neighbor edge structure) augmented by adding an  edge between any non-neighboring vertices~$x$ and~$y$ with probability
\begin{equation}
\label{E:1.1}
p_{x,y}:=1-\exp\{-\beta|x-y|^{-s}\}
\end{equation}
independently of all other edges. Here $\beta>0$ and~$s>0$ are parameters and $|\cdot|$ is any norm of choice. Our main point of interest is the behavior of the graph-theoretical distance~$\Ddis(x,y)$, defined as the minimal number of edges used in a path in that connects~$x$ to~$y$, in the limit as the separation between~$x$ and~$y$ tends to infinity. 

The question of distance scaling in long-range percolation has been studied quite intensely in the past and this has revealed five distinct regimes of typical behavior: $s<d$, $s=d$, $d<s<2d$, $s=2d$ and $s>2d$. Deferring the discussion of the specifics and references until the end of this section, let us focus directly on the regime $d<s<2d$. Here the first author~\cite{B1,B2} showed
\begin{equation}
\label{E:1.2}
\Ddis(0,x) = (\log |x|)^{\Delta+o(1)},\qquad |x|\to\infty,
\end{equation}
where
\begin{equation}
\label{E:1.3}
\Delta:=\frac1{\log_2(1/\gamma)}
\quad\text{for}\quad\gamma:=\frac s{2d}
\end{equation}
and where $o(1)\to0$ in probability. The proof worked for more general connection probabilities than \eqref{E:1.1}; in fact, it was enough to assume that~$p_{xy}=p_{0,x-y}=|x-y|^{-s+o(1)}$ as~$|x-y|\to\infty$.

The question we wish to resolve here is whether assuming the ``perfect'' scaling \eqref{E:1.1} yields a sharper version of the asymptotic \eqref{E:1.2}. Our first result in this regard is the subject of:

\begin{theorem}
\label{thm-1}
Consider the long-range percolation on~$\Z^d$ with connection probabilities \eqref{E:1.1} for~$\beta>0$ and~$s\in(d,2d)$ and let $\Ddis(x,y)$ denote the chemical distance between~$x$ and~$y$. There are~$c,C\in(0,\infty)$ depending only on~$\beta$,~$s$ and the underlying norm~$|\cdot|$ such that
\begin{equation}
\label{E:1.4}
\lim_{|x|\to\infty}P\bigl(c(\log |x|)^\Delta\le \Ddis(0,x)\le C(\log|x|)^\Delta\bigr)=1,
\end{equation}
where~$\Delta$ is as in \eqref{E:1.3}.
\end{theorem}

As soon as we accept \eqref{E:1.4}, a natural next step is the consideration of possible distributional limits of $\Ddis(0,x)/(\log|x|)^\Delta$ as~$|x|\to\infty$. We have been able to argue that if a distributional limit exists along a particular lattice direction, then it has to be non-random. Unfortunately, the proof of existence of the limit remains elusive, despite multiple attempts. Ultimately, this has led us to the consideration of a model on~$\R^d$ where progress can be made.

\begin{figure}[t]
\vglue-1mm
\centerline{\includegraphics[width=0.32\textwidth]{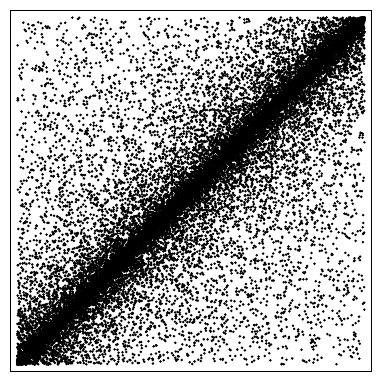}
\includegraphics[width=0.32\textwidth]{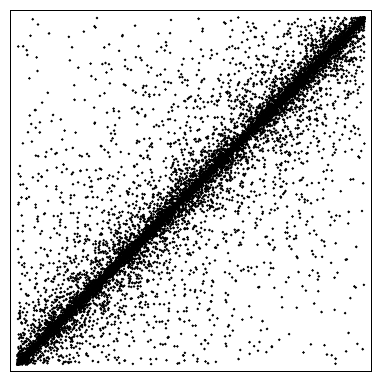}
\includegraphics[width=0.32\textwidth]{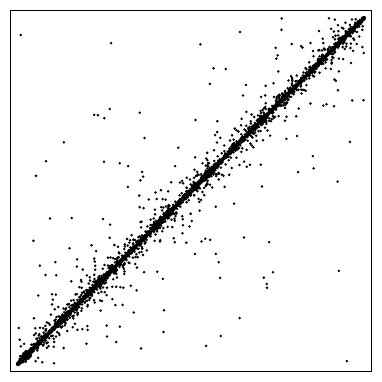}
}
\vglue0mm
\begin{quote}
\small 
\caption{
\label{fig-points}
\small
Sample plots of the symmetrized Poisson process $\Sym(\II_\beta)$ on~$\R$ for $\beta:=1$ and~$s:=1.1$ (left), $s:=1.4$ (middle) and~$s:=1.7$ (right). Only the restriction of the process to a box of side 5000 is shown.}
\normalsize
\end{quote}
\end{figure}

To define long-range percolation over~$\R^d$, fix~$\beta>0$ and consider a sample $\II_\beta$ from the Poisson process on $\R^d\times\R^d$ with ($\sigma$-finite) intensity measure
\begin{equation}
\label{E:1.15}
\mu_{s,\beta}(\textd x\,\textd y):=1_{\{\,|x|_2<|y|_2\}}\frac\beta{|x-y|^s}\textd x\,\textd y,
\end{equation}
where~$|\cdot|$ is the norm from \eqref{E:1.1} while~$|\cdot|_2$ is, here and henceforth, the Euclidean norm on~$\R^d$. Let us write~$\Sym(\II_\beta):=\II_\beta\cup\{(y,x)\colon (x,y)\in\II_\beta\}$ for the symmetrized version of~$\II_\beta$; see Fig.~\ref{fig-points}. We regard a ``point''~$(x,y)\in\Sym(\II_\beta)$ as an \textit{undirected edge} connecting~$x$ to~$y$. Given~$x,y\in\R^d$, we then proclaim
\begin{equation}
\label{E:1.6}
D(x,y):=\inf\biggl\{n+\sum_{i=0}^n|x_{i+1}-y_i|\colon n\ge0,\,\{(x_i,y_i)\colon i=1,\dots,n\}\subset\Sym(\II_\beta)\biggr\}
\end{equation}
with the convention $y_0:=x$ and~$x_{n+1}:=y$, to be the \textit{chemical distance} between points $x,y\in\R^d$ in the graph with edges~$\II_\beta$. 

\nopagebreak

We will at times refer to the sequence $\{(x_i,y_i)\colon i=1,\dots,n\}$ as a \textit{path} and call $x_{i+1}-y_i$ the $i$-th \textit{linear segment}.  Note that the infimum is over a non-empty set as the empty path, i.e., the one with~$n=0$ and one linear segment $y-x$, is always included. Note also that edges~$(x_i,y_i)$ with $|x_i-y_i|<1$ need not be considered as their removal decreases (thanks to the triangle inequality for the norm~$|\cdot|$ and the unit cost for each edge in the path) the expression in the infimum. In light of the local finiteness of~$\II_\beta$ away from $\{(x,y)\in\R^d\times\R^d\colon x=y\}$, this permits showing that the infimum is a.s.\ achieved. The main result of the present paper is then:

\begin{theorem}
\label{thm-2}
For each~$s\in(d,2d)$, each~$\beta>0$ (and $\gamma$ and~$\Delta$ as in \eqref{E:1.3}) and each choice of the norm~$|\cdot|$ there is a positive and continuous function~$\phi\colon(1,\infty)\to(0,\infty)$ satisfying
\begin{equation}
\label{E:1.7}
\phi(r^\gamma)=\phi(r),\quad r\ge1,
\end{equation}
such that for each~$x\in\R^d\smallsetminus\{0\}$,
\begin{equation}
\label{E:1.8}
\frac{D(0,rx)}{\phi(r)(\log r)^\Delta}\,\underset{r\to\infty}\longrightarrow\,1,\quad\text{in probability}.
\end{equation}
Moreover, $t\mapsto\phi(\texte^t)t^\Delta$ is convex throughout~$[0,\infty)$.
\end{theorem}

Note that we claim existence and continuity of~$\phi$ on~$(1,\infty)$ only. Actually,~$\phi$ admits a continuous extension to~$r=1$ if and only if it is constant.

\nopagebreak
\begin{figure}[t]
\vglue-1mm
\centerline{\includegraphics[width=0.5\textwidth]{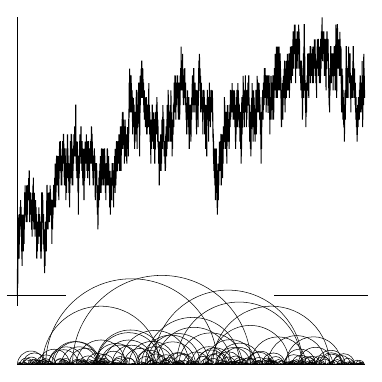}
}
\vglue0mm
\begin{quote}
\small 
\caption{
\label{fig-dist}
\small
Top curve: The distance from the origin in long-range percolation on~$\Z$ with $s:=1.8$ (and~$\beta:=1$) to points within (Euclidean) distance~10000. Arcs below: These depict the long edges in the corresponding portion of the graph. Note that the distance often dips down at points where a long edge lands.}
\normalsize
\end{quote}
\end{figure}

\subsection{Remarks and open questions}
We continue with some remarks and open questions. First we make:

\begin{observation}
The mode of convergence in \eqref{E:1.8} cannot be improved to almost sure.
\end{observation} 

\noindent
This is best seen in~$d=1$ by the following argument: A ball of radius~$r^\gamma$ centered at the origin will meet an edge of length of order~$r$ (in fact, even up to lengths~$r^{1/(2\gamma-1)}$) with a uniformly positive probability. If we parametrize the nearer endpoint of this edge as~$r^\gamma x$ and write~$ry$ for the farther endpoint of this edge, then the assumption of a.s.\ convergence in Theorem~\ref{thm-2} would tell us
\begin{equation}
\begin{aligned}
D(0,ry) &\le 1 + D(0,r^\gamma x)
\\
&=1+\phi(r^\gamma)(\log r^\gamma)^\Delta\bigl(1+o(1)\bigr)
\\
&=\frac12\phi(r)(\log r)^\Delta\bigl(1+o(1)\bigr),
\end{aligned}
\end{equation}
a contradiction with our very assumption. (We used \eqref{E:1.7} and the fact that~$\gamma^\Delta = \frac12$.) See Fig.~\ref{fig-dist}. 

Next, although this may not be quite apparent at first sight, the distance on~$\R^d$ is actually quite closely related to the chemical distance on~$\Z^d$. Indeed, replacing the Lebesgue measure on the right of~\eqref{E:1.15} by the counting measure on~$\Z^d$, the case when~$|\cdot|$ is the~$\ell^1$-norm on~$\R^d$ reduces \textit{exactly} to distance~$\Ddis(x,y)$ with~$x,y\in\Z^d$ connected with probability~$p_{x,y}$ as in~\eqref{E:1.1}. However, this does not seem to help in extending the sharp asymptotic \eqref{E:1.8} to the model on~$\Z^d$.

Another remark concerns the function~$\phi$ which encodes the dependence of the limit on~$\beta$ and the underlying norm~$|\cdot|$. We in fact believe:

\begin{conjecture}
The function~$\phi$ above is constant for each~$\beta>0$.
\end{conjecture}

\noindent
This is because~$\phi$ seems to appear largely as an artifact of our method which uses subadditivity arguments to relate the chemical distances at scales of the form $\{r^{\gamma^{-n}}\colon n\ge1\}$ for a fixed choice of~$r>1$. The growth rates of this sequence for two distinct~$r,r'\in[\texte^\gamma,\texte)$ are so incommensurate that the same proof would apply even if the intensity measure \eqref{E:1.15} were modulated depending on which of the two sequences~$|x-y|$ is closer to. In that situation, we would actually \textit{not} expect the corresponding~$\phi$ to take the same value at~$r$ and~$r'$. Unfortunately, we do not know how to turn this observation around to construct a proof of the above conjecture.

The dependence of $D(0,rx)$ on~$x$ is another interesting problem. As shown in Theorem~\ref{thm-2}, there is no such dependence in the leading order. Still, regardless on how the above conjecture gets resolved, formal expansions suggest:

\begin{conjecture}
For any~$x\ne0$ we have 
\begin{equation}
D(0,rx) =\phi(r)(\log r)^\Delta+\bigl(1+o(1)\bigr)\psi(r)(\log|x|)(\log r)^{\Delta-1}\,,
\end{equation}
where $o(1)\to0$ in probability as~$r\to\infty$ and where~$\psi$ is again a positive and continuous function satisfying the kind of ``periodicity'' requirement \eqref{E:1.7}.
\end{conjecture}

\noindent
This would in particular imply that balls in the chemical distance are close to those in the norm~$|\cdot|$. However, at this point we lack good ideas how to tackle this question rigorously.

Finally, our model admits an interesting generalization to a disordered setting obtained, on~$\Z^d$, by replacing \eqref{E:1.1} by
\begin{equation}
p_{xy}:=1-\exp\bigl\{-\beta W_xW_y|x-y|^{-s}\bigr\}
\end{equation}
for~$\{W_x\colon x\in\Z^d\}$ i.i.d.\ positive random variables. For this model, introduced by Deifen, van der Hofstad and Hooghiemstra~\cite{DvdHH} and with the scaling of distances studied by Deprez, Hazra and W\"uttrich~\cite{DHW}, different moment assumptions on the $W_x$'s imply different types of scalings of the chemical distance. It would be of interest to see whether any of the present proofs generalize to those situations as well.

\subsection{Earlier work and connections}
We will now give the promised connections to the existing literature on the scaling of the chemical distance in long-range percolation on~$\Z^d$. In the regime $s<d$ the chemical distance approaches a deterministic finite number at large spatial scales; namely, $\lceil\frac{d}{s-d}\rceil$ (Benjamini, Kesten, Peres and Schramm~\cite{BKPS}). When~$s=d$, the chemical distance between points at Euclidean distance~$N$ grows as~$(\log N)/\log\log N$ (Coppersmith, Gamarnik and Sviridenko~\cite{CGS}) while, as already mentioned, for~$d<s<2d$ we get \eqref{E:1.2} (Biskup~\cite{B1,B2}). For~$s>2d$, the chemical distance resumes linear scaling with the Euclidean distance (Berger~\cite{Berger-LRP}). As mentioned above, the asymptotics admit interesting generalizations to inhomogenous versions of long-range percolation (Deprez, Hazra and W\"uttrich~\cite{DHW}). 

The most interesting case is that of~$s=2d$, where the model is scale invariant. Some aspects of the~$d=1$ situation have been clarified already by Benjamini and Berger~\cite{Benjamini-Berger} but it was not until recently that Ding and Sly~\cite{DS} established the existence of an exponent~$\theta(\beta)\in(0,1)$ such that $D(0,N)\asymp N^{\theta(\beta)}$ in~$d=1$.  Interestingly, also here subadditivity arguments play a prominent role. The existence of a sharp asymptotic of the type established in Theorem~\ref{thm-2} remains open.

\subsection{Outline}
The rest of this note is organized as follows. In Section~\ref{sec2} we define the notion of a restricted distance~$\wt D$ and show (in Proposition~\ref{prop-2.2}) that it obeys a stochastic subadditivity bound that will drive all subsequent derivations in this paper. This bound produces~$\wt D$ at randomized locations and so its recursive use naturally leads, in Section~\ref{sec3}, to the consideration of a random variable~$W$ which is a fixed point of the randomization. Working with distances to multiples of~$W$ closes the recursion and permits extraction (in Proposition~\ref{prop-3.6}) of the limit asymptotic of~$r\mapsto\wt D(0,rW)$. A key technical step in this is the finiteness (derived in Lemma~\ref{lemma-2.6}) of the sum of conditional variances (given~$W$) of~$2^{-n}\wt D(0,r^{\gamma^{-n}}W)$ for any~$r\ge1$. Dini's theorem then yields uniformity of such estimates in~$r$ which in turn allows for elimination of~$W$ from the argument of the distance. In Section~\ref{sec4} we then show that the same asymptotic applies to distance~$D$ as well.

\section{Restricted distance}
\label{sec2}\nopagebreak\noindent
We are now ready to commence the proofs. A majority of the work will be done directly for the model on~$\R^d$ although we do use the model on~$\Z^d$ in the proof of positivity of~$\phi$. In this section we focus on an auxiliary quantity, called the restricted distance, that is better behaved under subadditivity arguments. We will return to the full distance in Section~\ref{sec4}.

\subsection{Definition and comparisons}
Let us write~$B(x,r):=\{y\in\R^d\colon |x-y|<r\}$ for the open ball in the norm~$|\cdot|$. Given~$x,y\in\R^d$ we then define their \emph{restricted distance} by constraining the infimum in \eqref{E:1.6} to paths that do not leave the ball~$B(x,2|x-y|)$, i.e.,
\begin{equation}
\label{E:2.1}
\wt D(x,y):=\inf\Biggl\{n+\sum_{i=0}^n|x_{i+1}-y_i|\,\colon\,
\begin{aligned}
&n\ge0,\,\{(x_i,y_i)\colon i=1,\dots,n\}\subset\II_\beta
\\
&x_i,y_i\in B\bigl(x,2|x-y|\bigr)\,\,\forall i=1,\dots,n
\end{aligned}
\Biggr\}\,,
\end{equation}
where, as before, we set~$y_0:=x$ and~$x_{n+1}:=y$. We caution the reader that~$\wt D(\cdot,\cdot)$ is not a metric as it is neither symmetric nor obeying the triangle inequality. The following properties of the restricted distance will be important in the sequel:

\begin{lemma}
\label{lemma-2.1}
Let~$D(x,y)$ be as in \eqref{E:1.6} and $\wt D(x,y)$ as in \eqref{E:2.1}. Then
\settowidth{\leftmargini}{(1111)}
\begin{enumerate}
\item[(1)] for any $x,y\in\R^d$,
\begin{equation}
\label{E:2.2}
D(x,y)\le\wt D(x,y)\le|x-y|,
\end{equation}
\item[(2)] the law of $\wt D$ is translation invariant,
\begin{equation}
\label{E:2.3}
\bigl\{\wt D(x,y)\colon x,y\in\R^d\bigr\}\laweq \bigl\{\wt D(x+z,y+z)\colon x,y\in\R^d\bigr\},\qquad z\in\R^d,
\end{equation}
\item[(3)] $x\mapsto\wt D(0,x)$ is stochastically continuous in~$x$ (i.e., the law of $\wt D(0,x)$ is continuous  in~$x$ in the topology of weak convergence of measures), and
\item[(4)] for any~$x,y,\tilde x,\tilde y\in\R^d$
\begin{equation}
\label{E:2.4}
|x-\tilde x|>2|x-y|+2|\tilde x-\tilde y|
\quad\Rightarrow\quad
\wt D(x,y)\independent\wt D(\tilde x,\tilde y)
\end{equation}
\end{enumerate}
\end{lemma}
 
\begin{proofsect}{Proof}
The inequalities \eqref{E:2.2} are checked by comparison of \eqref{E:1.6} with \eqref{E:2.1} and the fact that the path with no edges is included on the right of \eqref{E:2.1}. The translation invariance in~(2) is a consequence of the corresponding property of the intensity measure \eqref{E:1.15}. 

In order to prove~(3), consider a path minimizing~$\wt D(0,x)$. (Such a path exists as $B(0,2|x|)$ contains only a finite number of edges of length in excess of one, a.s.) The continuity of the law of the underlying point process and the fact that $B(0,2|x|)$ is open ensure that the minimizing path is a.s.\ unique and that the same sequence of edges are used by the minimizer of~$\wt D(0,x+z)$ for all~$|z|$ sufficiently small. It follows that, for every~$y\in\R^d$, the map~$x\mapsto\wt D(0,x)$ is continuous at~$y$ a.s. This yields stochastic continuity via the Bounded Convergence Theorem.

The independence claimed in \eqref{E:2.4} follows from the independence of Poisson processes over disjoint sets.
\end{proofsect}

Let us write~$D_\beta$ if need arises to mark explicitly the dependence of the law random variable~$D$ on~$\beta$. The following comparisons then hold:

\begin{lemma}
\label{lemma-4.3}
For all~$\beta>0$, all~$a\ge1$ and all~$x\in\R^d$,
\begin{equation}
\label{E:1.7w}
D_\beta(0,ax) \,\,\overset{\text{\rm law}}\le\,\,
D_{a^{s-2d}\beta}(0,ax)
\,\,\overset{\text{\rm law}}\le\,\,
a D_\beta(0,x).
\end{equation}
The same conclusions apply to the restricted distance~$\wt D_\beta$ as well.
\end{lemma}

\begin{proofsect}{Proof}
Let $\II_\beta=\{(x_i,y_i)\colon i\in\N\}$ denote a sample from the point process on~$\R^d\times\R^d$ with intensity measure \eqref{E:1.15}. For any~$a>0$, the process $\II'_\beta:=\{(ax_i,ay_i)\colon i\in\N\}$ is then equidistributed to~$\II_{\beta(a)}$ where $\beta(a):=a^{s-2d}\beta$. Pick a path $\pi$ connecting $0$ to~$x$ using the edges in~$\Sym(\II_\beta)$ and let~$n(\pi)$ denote the number of edges and~$\ell(\pi)$ the total length of the linear segments in~$\pi$. Now consider the path~$\pi'$ built using the corresponding edges in~$\Sym(\II_\beta')$, and let~$n(\pi')$ and~$\ell(\pi')$ denote the corresponding quantities for~$\pi'$. Then
\begin{equation}
n(\pi')=n(\pi)\quad\text{and}\quad \ell(\pi')=a \ell(\pi).
\end{equation}
Assuming~$a\ge1$, it follows that
\begin{equation}
\label{E:2.7a}
n(\pi')+\ell(\pi')\le a\bigl[n(\pi)+\ell(\pi)\bigr].
\end{equation}
The left-hand side is bounded by $D(0,ax)$ from below; optimizing over~$\pi$ then implies the inequality on the right of \eqref{E:1.7w}. The left inequality in \eqref{E:1.7w} is a consequence of Poisson thinning: As $\beta(a)\le\beta$ for $a\ge1$, the process with parameter~$\beta(a)$ can be realized as a pointwise subset of the process with parameter~$\beta$. In this coupling, every path contributing to $\wt D_{\beta(a)}(0,ax)$ will contribute to~$\wt D_\beta(0,ax)$ as well. 

The inequalities extend to~$\wt D$ as the additional restriction imposed on paths there scales proportionally to the distance between endpoints.
\end{proofsect}

For distance~$D$ we can also get comparisons under rotations:

\begin{lemma}
\label{lemma-4.3a}
For each~$\epsilon>0$ there is~$\delta>0$ such that for all~$x,y\in\R^d\smallsetminus\{0\}$,
\begin{equation}
\label{E:2.6ab}
|x|_2=|y|_2\quad\&\quad\frac{|x-y|_2}{|x|_2}<\delta\quad\Rightarrow\quad
D_{(1+\epsilon)\beta}(0,x)\,\,\overset{\text{\rm law}}\le\,\,(1+\epsilon)D_\beta(0,y).
\end{equation}
\end{lemma}

\begin{proofsect}{Proof}
Thanks to all norms on~$\R^d$ being continuous with respect to one another, for each~$\epsilon>0$ there is~$\delta>0$ such that for any rotation~$R\in\text{SO}(d)$ which is close to the identity in the sense that~$|Rx-x|_2<\delta|x|_2$ for all non-zero~$x\in\R^d$,  we have
\begin{equation}
\label{E:2.9a}
(1+\epsilon)|x|\ge\vert Rx\vert\ge (1+\epsilon)^{-1/s}\vert x\vert,\qquad x\in\R^d\smallsetminus\{0\}.
\end{equation}
The inequality on the right shows that $\mu_{s,\beta(1+\epsilon)}-\mu_{s,\beta}\circ R^{-1}$ is a positive measure. The additivity of Poisson processes implies that a sample~$\II_{\beta(1+\epsilon)}$ from the Poisson process with the intensity~$\mu_{s,\beta(1+\epsilon)}$ contains a sample $\II_\beta'$ from the process with intensity~$\mu_{s,\beta}$ rotated by~$R$. Pick a path~$\pi'$ using the edges in~$\Sym(\II_\beta')$ from~$0$ to~$R^{-1}x$ and let~$\pi$ be its rotation by~$R$. Then, in the notation from the previous proof,~$n(\pi)=n(\pi')$ while, by the left inequality in \eqref{E:2.9a}, $\ell(\pi)\le(1+\epsilon)\ell(\pi')$. Optimizing over~$\pi'$ we get
\begin{equation}
D_{\beta(1+\epsilon)}(0,x)\,\,\overset{\text{\rm law}}\le\,\, (1+\epsilon) D_\beta(0,R^{-1}x)
\end{equation}
for every~$x\in\R^d$.
Realizing~$y$ as~$R^{-1}x$, this yields \eqref{E:2.6ab}. 
\end{proofsect}

We can even get comparisons with the distance on~$\Z^d$, writing again $\Ddis_\beta$ to denote the distance on~$\Z^d$ with connection probabilities \eqref{E:1.1} for parameter~$\beta$:

\begin{lemma}
\label{lemma-3th}
For each~$\beta>0$ there is~$c=c(\beta)\in(0,1]$ such that for all~$x,y\in\Z^d$,
\begin{equation}
c\Ddis_{c^{-1}\beta}(x,y) \,\,\overset{\text{\rm law}}\le\,\, D_\beta(x,y) \,.
\end{equation}
\end{lemma}

\begin{proofsect}{Proof}
Denote~$B:=[-1/2,1/2)^d$ and define the coupling of the process on~$\Z^d$ and the process on~$\R^d$ as follows. Given a sample~$\II_\beta$ from the Poisson process with intensity~$\mu_{s,\beta}$, place an edge between distinct non-nearest neighbors~$x\in\Z^d$ and~$y\in\Z^d$ whenever there is an edge $(x',y')\in\Sym(\II_\beta)$ with $x'-x,y'-y\in B$. Distinct vertices~$x,y\in\Z^d$ are then connected by an edge with probability 
\begin{equation}
1-\exp\Bigl\{-\beta\int_{B\times B}\frac{\textd z\,\textd z'}{|x-y+z-z'|^{s}}\Bigr\}
\end{equation}
independent of all other edges. Note that (since~$s>d$) the integral diverges for any two~$x,y\in\Z^d$ within~$\ell^\infty$-distance one which ensures these points are connected almost surely. As is readily checked, the resulting process on~$\Z^d$ stochastically dominates the process defined in \eqref{E:1.1} with~$\beta$ multiplied by a sufficiently large constant.

Now consider a path~$\pi$ contributing to~$D_\beta(x,y)$ and use the above coupling to project it to a path~$\pi'$ on~$\Z^d$ while replacing each linear segments of~$\pi$ by a shortest nearest-neighbor path on~$\Z^d$ between the corresponding vertices on~$\Z^d$. An edge in~$\pi$ then gives rise to an edge in~$\pi'$ or no edge at all. A linear segment in~$\pi$ of length~$L$ corresponds to a ``segment'' on~$\Z^d$ of $\ell^1$-distance~$L'$ between the endpoints or no segment at all. The fact that the $\ell^1$-distance is comparable with the norm~$|\cdot|$ ensures that $L\ge cL'$ for some~$c>0$ small enough. The claim then follows.
\end{proofsect}

We will find the lower bound by distance on~$\Z^d$ particularly useful in light of the following result by the first author that itself draws on earlier work by Trapman~\cite{Trapman}:

\begin{theorem}
\label{thm-4th}
For each~$\beta>0$ there are~$c_1,c_2\in(0,\infty)$ such that for all~$n\ge1$ and all~$x\in\Z^d$,
\begin{equation}
\label{E:2.9th}
P\bigl(\Ddis(0,x)\le n\bigr)\le c_1\frac{\texte^{c_2n^{1/\Delta}}}{|x|^s}\,.
\end{equation}
\end{theorem}

\begin{proofsect}{Proof}
This is proved by following, nearly verbatim, the proof of~\cite[Theorem~3.1]{B2} while setting~$s':=s$ and~$\Delta':=\Delta$. (Note that~$s'$ is introduced in~\cite{B2} in order to reduce the asymptotic form $|x-y|^{-s+o(1)}$ assumed there for~$p_{xy}$ to the sharp asymptotic~\eqref{E:1.1} with~$s'$ instead of~$s$. The rest of the proof then uses the sharp asymptotic form of~$p_{xy}$ only.)
\end{proofsect}

From here we get one half of Theorem~\ref{thm-1} of the present paper:

\begin{corollary}
\label{cor-2.5}
For each~$\beta>0$ there is~$c=c(\beta)>0$ such that
\begin{equation}
\lim_{|x|\to\infty}\,P\bigl(\Ddis(0,x)\le c(\log|x|)^\Delta\bigr)=0.
\end{equation}
\end{corollary}

\begin{proofsect}{Proof}
Substitute $n:=c(\log|x|)^\Delta$ into \eqref{E:2.9th} and observe that, thanks to~$s>d$, the resulting probability is summable on~$x\in\Z^d$ once~$c$ is sufficiently small. This implies the claim.
\end{proofsect}

\subsection{Subadditivity bound}
Our next task is to derive a subadditivity relation for the restricted distance. This relation will play a fundamental role in all derivations to come. We remark that~$D$, being a metric, satisfies the ``ordinary'' subadditivity estimate
\begin{equation}
D\bigl(0,(n+m)x\bigr)\le D(0,nx)+D\bigl(nx,(n+m)x\bigr).
\end{equation}
However, this estimate is not useful for our purposes because $n\mapsto D(0,nx)$ turns out to be sublinear a.s. Our subadditivity bound will thus have to be tailored to the polylogarithmic growth of~$x\mapsto D(0,x)$. It will also be derived only for the restricted distance because that, unlike~$D$, obeys the independence statement in Lemma~\ref{lemma-2.1}(4).

\begin{proposition}[Subadditivity for restricted distance]
\label{prop-2.2}
Fix $\eta\in(0,1)$ and let~$Z,Z'$ be i.i.d.\ $\R^d$-valued random variables with common law
\begin{equation}
\label{E:2.6a}
P(Z\in B)=\sqrt{\eta\beta}\int_B\texte^{-\eta\beta c_0|z|^{2d}}\textd z\,,
\end{equation}
where
\begin{equation}
\label{E:2.5}
c_0:=\int_{|z|^{2d}+|\tilde z|^{2d}\le1}\textd z\,\textd\tilde z\,.
\end{equation}
Let~$\wt D'$ be an independent copy of~$\wt D$ with both quantities independent of~$Z$ and~$Z'$.
For each $\gamma_1,\gamma_2\in(0,\frac12(1+\gamma))$ with $\gamma_1+\gamma_2=2\gamma$, there are $c_1,c_2\in(0,\infty)$ and, for each~$x\in\R^d$, there is an event~$A(x)\in\sigma(Z,Z')$ with
\begin{equation}
\label{E:2.6}
P\bigl(A(x)\bigr)\le c_1\texte^{-c_2|x|^{\theta}}
\end{equation}
for $\theta:=2d[\frac{1+\gamma}2-\gamma_1\vee\gamma_2]$ such that
\begin{equation}
\label{E:2.7}
\wt D(0,x)\overset{\text{\rm law}}\le\, \wt D\bigl(0,|x|^{\gamma_1} Z\bigr)+\wt D'\bigl(0,|x|^{\gamma_2} Z'\bigr)+1+|x|1_{A(x)}
\end{equation}
holds true for every~$x\in\R^d$.
\end{proposition}

\begin{remark}
It may not appear obvious that the above choice of~$c_0$ makes \eqref{E:2.6a} a probability; this will be seen from formula \eqref{E:2.13} in the proof below. We will use this proposition mostly in the case when~$\gamma_1=\gamma_2=\gamma$. The main reason for our consideration of the more general setting is the proof of continuity of the limit in Theorem~\ref{thm-2} which requires (small but non-trivial) perturbations about the symmetric case as well. The choice of~$\eta$ will be immaterial in what follows. We will therefore suppress~$\eta$ from the notation wherever possible.
\end{remark}

\begin{proofsect}{Proof of Proposition~\ref{prop-2.2}} 
The main idea of the proof is simple: We first pick an edge $(X,Y)$, with~$X$ closest to~$0$ and~$Y$ closest to~$x$ according to criteria to be specified later. Then we pick a shortest path from~$0$ to~$X$ and a shortest path from~$x$ to~$Y$, demanding in addition that the first path stay in~$2|X|^{\gamma_1}$-neighborhood of~$0$ and the second in $2|x-Y|^{\gamma_2}$-neighborhood of~$x$. Assuming $|x|\gg1$, concatenating the two paths with~$(X,Y)$ yields a path from~$0$ to~$x$ not leaving $2|x|$-neighborhood of~$0$. Writing~$Z$ for~$|x|^{-\gamma_1} X$ and~$Z'$ for~$|x|^{-\gamma_2}(Y-x)$, a pointwise version of \eqref{E:2.7} follows. A key technical point is to choose the selection criteria for~$(X,Y)$ to ensure independence of~$Z$ and~$Z'$ and (conditionally on~$X$ and $Y$) the distances $\wt D(0,X)$ and $\wt D(x,Y)$.

\smallskip
Fix~$\eta\in(0,1)$. There is nothing to prove when~$x=0$ so let us also assume that~$x\in\R^d\smallsetminus\{0\}$. Recall that~$a\vee b$ denotes~$\max\{a,b\}$. The proof comes in three steps.

\medskip
\noindent
\textsl{STEP 1: Construction of~$(X,Y)$}: We begin by constructing the aforementioned edge.
Note that, for any~$\tilde x,\tilde y\in\R^d$ with $|\tilde x|\vee|\tilde y-x|\le|x|^{\frac12(1+\gamma)}$ we have
\begin{equation}
\label{E:2.8}
|\tilde x-\tilde y|\le|x+\tilde x+\tilde y-x|
\le|x|+2|x|^{\frac12(1+\gamma)}=|x|\bigl(1+2|x|^{-\frac12(1-\gamma)}\bigr).
\end{equation}
Recalling the intensity measure~$\mu_{s,\beta}$ from \eqref{E:1.15}, define
\begin{equation}
\label{E:2.9}
\mu_{s,\beta}'(\textd\tilde x\,\textd\tilde y):=\eta\beta\1_{\{|\tilde x|_2<|\tilde y|_2\}}
\1_{\bigl\{|\tilde x|\vee |\tilde y-x|\le |x|^{\frac12(1+\gamma)}\bigr\}}
\frac{\textd\tilde x\,\textd\tilde y}{|x|^s}.
\end{equation}
Then, as soon as~$x$ is so large that $1+2|x|^{-\frac12(1-\gamma)}\le\eta^{-1/s}$ (recall that~$\gamma\in(0,1)$ and~$\eta<1$), the inequality~\eqref{E:2.8} ensures that $\mu_{s,\beta}'':=\mu_{s,\beta}-\mu_{s,\beta}'$ is a positive measure. This permits us to represent~$\II_\beta$ as the sum of two independent Poisson processes~$\II_\beta'$ and~$\II_\beta''$ with intensities $\mu_{s,\beta}'$ and~$\mu_{s,\beta}''$, respectively. Considering also the measure
\begin{equation}
\label{E:2.10}
\mu_{s,\beta}'''(\textd\tilde x\,\textd\tilde y):=\eta\beta\Bigl(1-\1_{\{|\tilde x|_2<|\tilde y|_2\}}\1_{\bigl\{|\tilde x|\vee |\tilde y-x|\le |x|^{\frac12(1+\gamma)}\bigr\}}\Bigr)
\frac{\textd\tilde x\,\textd\tilde y}{|x|^s},
\end{equation}
let~$\II_\beta'''$ denote a sample of the Poisson process with intensity~$\mu_{s,\beta}'''$. We regard $\II_\beta'$, $\II_\beta''$ and $\II_\beta'''$ as independent of one another.

As is directly checked from \twoeqref{E:2.9}{E:2.10}, $\II_\beta'\cup\II_{\beta}'''$ is a homogeneous Poisson process on~$\R^d\times\R^d$ with density~$\eta\beta|x|^{-s}$ and so, in particular, $\II_\beta'\cup\II_{\beta}'''\ne\emptyset$ a.s. The process is also locally finite and so there is (a.s.) a unique pair  $(X,Y)\in\II_\beta'\cup\II_{\beta}'''$ minimizing the function
\begin{equation}
\label{E:2.11}
f_x(\tilde x,\tilde y):=\bigl(|x|^{-\gamma_1}|\tilde x|\bigr)^{2d}+\bigl(|x|^{-\gamma_2}|\tilde y-x|\bigr)^{2d}.
\end{equation}
Setting
\begin{equation}
\label{E:2.12}
Z_x:=|x|^{-\gamma_1}X\quad\text{and}\quad Z_x':=|x|^{-\gamma_2}(Y-x),
\end{equation}
and noting that $d\gamma_1+d\gamma_2=s$, the law of $(Z_x,Z_x')$ is given by
\begin{equation}
\label{E:2.13}
P\bigl(Z_x\in \textd z,\,Z_x'\in\textd z'\bigr)
=\eta\beta\exp\Bigl\{-\eta\beta\int \textd \tilde z\,\textd\tilde z'\1_{\{|\tilde z|^{2d}+|\tilde z'|^{2d}\le|z|^{2d}+|z'|^{2d}\}}\Bigr\}
\textd z\,\textd z'\,.
\end{equation}
Scaling the variables in the inner integral by $(|z|^{2d}+|z'|^{2d})^{\frac1{2d}}$ and invoking \eqref{E:2.5} shows that $(Z_x,Z_x')$ are i.i.d.\ with law as in~\eqref{E:2.6a}.

\medskip
\noindent
\textsl{STEP 2: Definition of~$A(x)$ and pointwise inequality}:
We will now define~$A(x)$ and prove a pointwise version of the inequality \eqref{E:2.7}. For~$x$ so large that $1+2|x|^{-\frac12(1-\gamma)}\le\eta^{-1/s}$ and $4|x|^{\frac12(1+\gamma)}<|x|$ hold true, we set
\begin{equation}
A(x):=\bigl\{|Z_x|> |x|^{\frac12(1+\gamma)-\gamma_1}\bigr\}\cup\bigl\{|Z_x'|> |x|^{\frac12(1+\gamma)-\gamma_2}\bigr\}
\end{equation} 
and otherwise set $A(x)$ to be the sample space carrying the three Poisson processes above. We now claim the pointwise inequality
\begin{equation}
\label{E:2.15}
\wt D(0,x)\le\wt D\bigl(0,|x|^{\gamma_1} Z_x\bigr)+\wt D\bigl(x,x+|x|^{\gamma_2} Z_x'\bigr)+1+|x|1_{A(x)},
\end{equation}
where all instances of~$\wt D$ are defined using the symmetrized version of~$\II_\beta=\II_\beta'\cup\II_\beta''$.
Since $\wt D(0,x)\le|x|$, \eqref{E:2.15} holds whenever~$A(x)$ occurs and so we just need to verify \eqref{E:2.15} under the conditions
\begin{equation}
\label{E:2.16}
1+2|x|^{-\frac12(1-\gamma)}\le\eta^{-1/s},\quad 
4|x|^{\frac12(1+\gamma)}<|x|\quad
\text{and}\quad |X|,|x-Y|\le |x|^{\frac12(1+\gamma)}\,.
\end{equation}
Noting that $\mu_{s,\beta}'$ and~$\mu_{s,\beta}'''$ have disjoint supports, the last two conditions ensure $(X,Y)\in\II_\beta'\subset\II_\beta$ a.s.\ and so~$(X,Y)$ is allowed to enter a path contributing to the distance on the left of \eqref{E:2.15}. Fix any~$\epsilon>0$ and consider a path in $B(0,2|X|)$ from~$0$ to~$X$ of length at most~$\wt D(0,X)+\epsilon$ and then a path in $B(x,2|x-Y|)$ from~$x$ to~$Y$ of length at most~$\wt D(x,Y)+\epsilon$. Since \eqref{E:2.16} ensures
\begin{equation}
B(0,2|X|)\subseteq B(0,2|x|^{\frac12(1+\gamma)})\subseteq B(0,2|x|)
\end{equation}
and
\begin{equation}
B(x,2|x-Y|)\subseteq B(x,2|x|^{\frac12(1+\gamma)})\subseteq x+B(0,|x|)\subseteq B(0,2|x|),
\end{equation}
concatenating the former path with edge $(X,Y)$ and then adjoining the latter path after~$Y$, we get a path  contributing potentially to the infimum defining~$\wt D(0,x)$ and having length at most $\wt D(0,X)+\wt D(x,Y)+2\epsilon+1$. As~$\epsilon$ was arbitrary, \eqref{E:2.15} follows via \eqref{E:2.12}.

\medskip
\noindent
\textsl{STEP 3: Reduction to independent variables}: 
Let us now see how \eqref{E:2.15} reduces to \eqref{E:2.7}. Enlarge the probability space so that it holds two independent copies~$\wt D'$ of~$\wt D''$ of random variable~$\wt D$, that are independent of the processes~$\II'_\beta,\II''_\beta,\II'''_\beta$ and thus of the random objects~$\wt D$, $Z_x$ and~$Z_x'$. Under the restrictions on~$x$ from~\eqref{E:2.16} we have
\begin{equation}
B\bigl(0,2|x|^{\frac12(1+\gamma)}\bigr)\cap B\bigl(x,2|x|^{\frac12(1+\gamma)}\bigr)=\emptyset.
\end{equation}
It follows that, conditional on $A(x)^\cc$, the triplet of families of random variables
\begin{equation}
\label{E:2.29}
\bigl\{\wt D(0,z)\colon |z|< 2|x|^{\frac12(1+\gamma)}\bigr\},\quad 
\bigl\{\wt D(x,x+z)\colon |z|< 2|x|^{\frac12(1+\gamma)}\bigr\}\quad\text{and}\quad \{X,Y\}
\end{equation}
are independent. Moreover, $\wt D(x,x+z) \laweq \wt D''(0,z)$ by translation symmetry of the underlying process. Since, as before, \eqref{E:2.15} holds trivially when~$A(x)$ occurs, it suffices to check \eqref{E:2.7} conditionally on~$A(x)^\cc$. In that case the independence of the objects in \eqref{E:2.29} permits us to swap $\wt D(0,|x|^{\gamma_1}Z_x)$ for $\wt D'(0,|x|^{\gamma_1}Z_x)$ and $\wt D(x,x+|x|^{\gamma_2} Z_x')$ for~$\wt D''(0,|x|^{\gamma_2} Z_x')$ without affecting the (conditional) law of the right-hand side of \eqref{E:2.15}. Then \eqref{E:2.7} follows from \eqref{E:2.15}.

In order to complete the proof, it remains to verify the bound \eqref{E:2.6}. Assuming the first two conditions in \eqref{E:2.16} hold, $A(x)$ will occur only if one of~$Z_x$ or~$Z_x'$ exceeds the stated bounds. The formula \eqref{E:2.13} then readily shows \eqref{E:2.6} in this case. We then adjust the constant~$c_1$ so that \eqref{E:2.6} holds even when the first two conditions in \eqref{E:2.16} fail.
\end{proofsect}

\section{Limit considerations}
\label{sec3}\nopagebreak\noindent
The main goal of this section is to establish the limit claim from Theorem~\ref{thm-1} for the restricted distance.
Due to our lack of a suitable substitute for the Subadditive Ergodic Theorem, we will extract the result by controlling the expectation and the variance of the restricted distance. Throughout this section we fix~$\beta>0$ and~$\eta\in(0,\infty)$ and suppress them from all formal statements. 

\subsection{Convergence along doubly-exponential sequences} 
We begin by noting that iterations of \eqref{E:2.7} naturally lead to the consideration of  \emph{randomized} locations to which the restricted distance is to be computed. In order to get a closed-form expression, a natural idea is to work with a fixed point of the randomization. This leads to:

\begin{lemma}
\label{lemma-3.1}
Let $Z_0,Z_1,\dots$ be i.i.d.\ copies of the random variable from \eqref{E:2.6a}. Then the infinite product in
\begin{equation}
\label{E:2.19}
W:=Z_0\prod_{k=1}^\infty |Z_k|^{\gamma^k},
\end{equation}
converges in~$(0,\infty)$ a.s. Moreover,~$W$ has continuous, non-vanishing probability density and has all moments. Furthermore, if~$Z$ has the law as in \eqref{E:2.6a}, then\begin{equation}
\label{E:3.2ua}
Z\independent W\quad\Rightarrow\quad |W|^\gamma Z\laweq W
\end{equation}
\end{lemma}

\begin{proofsect}{Proof}
The random variable~$|\log Z|$ has exponential tails and so $k\mapsto |\log Z_k|$ grows at most polylogarithmically fast a.s. Since~$k\mapsto\gamma^k$ decays exponentially, the infinite product converges to a number in~$(0,\infty)$ a.s. This, along with the fact that~$Z$ has continuous and positive density, implies that~$W$ has continuous and positive density as well. 

To control the upper tail of~$W$, observe that  $\sum_{k\ge0}(k+1)\gamma^k=(1-\gamma)^{-2}$. Hence, if~$|W|> t^{(1-\gamma)^{-2}}$, then we must have $|Z_k|> t^{k+1}$ for at least one~$k\ge0$. Hereby we get
\begin{equation}
P\bigl(|W|>t\bigr)\le\sum_{k\ge0}P\bigl(|Z_k|>t^{(k+1)(1-\gamma)^2}\bigr),\qquad t>0.
\end{equation}
As the tails of~$Z$ are no heavier than Gaussian in all~$d\ge1$, the claim follows. The distributional identity $|W|^\gamma Z\laweq W$ for~$Z\independent W$ is checked directly from the definition of~$W$.
\end{proofsect}

The identity \eqref{E:3.2ua} shows that~$W$ is indeed a fixed point for the \textit{random} arguments of the restricted distance under iterations of \eqref{E:2.7}. This enables us to prove our first limit claim:

\begin{lemma}
\label{lemma-3.2}
Assume~$W$ from \eqref{E:2.19} is independent of~$\wt D$. Then for each~$r\ge1$, the limit
\begin{equation}
\label{E:2.20}
L(r):=\lim_{n\to\infty}\frac{E \,\wt D(0,r^{\gamma^{-n}}W)}{2^n},
\end{equation}
exists. Moreover, $r\mapsto L(r)$ is upper-semicontinuous on~$[1,\infty)$ and positive on~$(1,\infty)$.
\end{lemma}

\begin{proofsect}{Proof}
We will apply Proposition~\ref{prop-2.2} for the choices~$\gamma_1=\gamma_2=\gamma$.
Let~$W$ be the random variable independent of~$\wt D$, $\wt D'$, $Z$ and~$Z'$ in \eqref{E:2.7}. Plugging~$r^{\gamma^{-n}}W$ for~$x$ in \eqref{E:2.7} and invoking $|W|^\gamma Z\laweq W$ along with the bound \eqref{E:2.6} yields
\begin{equation}
\label{E:3.4th}
E \,\wt D(0,r^{\gamma^{-n}}W)\le 2E \,\wt D(0,r^{\gamma^{-n+1}}W) + c,
\end{equation}
where $c:=1+c_1\sup_{x\in\R^d}|x|\texte^{-c_2|x|^\theta}$, for~$c_1$, $c_2$ and~$\theta$ as in Proposition~\ref{prop-2.2}.
This shows that
\begin{equation}
a_n(r):=2^{-n}\bigl[E \,\wt D(0,r^{\gamma^{-n}}W)+c\bigr]
\end{equation}
 is non-increasing and, being non-negative,~$\lim_{n\to\infty}a_n(r)$ exists. The limit in \eqref{E:2.20} then exists as well and takes the same value. By Lemma~\ref{lemma-2.1}(3), $r\mapsto a_n(r)$ is continuous and so and so, being a decreasing limit of continuous functions, $r\mapsto L(r)$ is upper semicontinuous. The positivity of~$L(r)$ for~$r>1$ follows from Lemma~\ref{lemma-3th}, Theorem~\ref{thm-4th} and the fact that $[\log(r^{\gamma^{-n}})]^\Delta = 2^{n}(\log r)^\Delta$.
\end{proofsect}

We now augment the convergence of expectations to:

\begin{proposition}
\label{prop-2.5}
For any~$r\ge1$ and Lebesgue a.e.\ $x\in\R^d$, 
\begin{equation}
\label{E:3.6th}
\frac{\wt D(0,r^{\gamma^{-n}}x)}{2^n}\,\,\underset{n\to\infty}\longrightarrow\,\, L(r)\quad\text{\rm$P$-a.s.}
\end{equation}
In particular, $L(r)$ defined in \eqref{E:2.20} does not depend on the choice of~$\eta$.
\end{proposition}

The main ingredient of the proof is:

\begin{lemma}
\label{lemma-2.6}
Suppose~$\wt D$, $Z$ and~$W$ are independent with distribution as above.
Let~$\sigma(W)$ denote the sigma algebra generated by~$W$. Then for any~$r\ge1$,
\begin{equation}
\label{E:1.14th}
\sum_{n=1}^\infty E\biggl(\Var\Bigl(2^{-n}\wt D\bigl(0,r^{\gamma^{-n}}Z|W|^\gamma\bigr)\,\Big|\,\sigma(W)\Bigr)\biggr)<\infty\,.
\end{equation}
\end{lemma}

\begin{proofsect}{Proof}
Fix~$r\ge1$. Plugging~$x:=rW$ in \eqref{E:2.7}, squaring both sides and taking expectations we get
\begin{equation}
\label{E:3.9a}
E\bigl(\wt D(0,rW)^2\bigr)\le 2 E\bigl(\,\wt D(0,r^\gamma W)^2\bigr)+2E\Bigl(\bigl[E\bigl(\wt D(0,r^\gamma Z|W|^\gamma)\,\big|\,\sigma(W)\bigr)\bigr]^2\Bigr)+F_0\,,
\end{equation}
where, using $\wt D(0,x)\le|x|$ and $Z|W|^\gamma\laweq W$, the error term is given by
\begin{equation}
F_0=F_0(r):=E\bigl([1+r|W|1_{A(rW)}]^2\bigr)+4E\wt D(0, r^{\gamma}W)+4E\bigl[(r|W|)^{1+\gamma}1_{A(rW)}\bigr]\,.
\end{equation}
Next we rewrite the second term on the right of \eqref{E:3.9a} using conditional variance, and then subtract suitable terms on both sides to get
\begin{multline}
\quad
\Var(\wt D(0,rW))\le 2\Var(\wt D(0,r^\gamma W))+2\Var\Bigl(\bigl[E(\wt D(0,r^\gamma Z|W|^\gamma)|\sigma(W))\bigr]^2\Bigr)
\\
+4E\bigl(\wt D(0,r^\gamma W)\bigr)^2-E\bigl(\wt D(0,rW)\bigr)^2+F_0.
\quad
\end{multline}
Replacing~$W$ by $Z|W|^\gamma$ in the first two variances above and using the standard identity
\begin{equation}
\Var(X)=E\bigl(\Var(X|Y)\bigr)+\Var\bigl(E(X|Y)\bigr)
\end{equation} 
yields
\begin{multline}
\label{E:3.13}
\quad
\Var\bigl(E(\wt D(0,rZ|W|^\gamma)\,\big|\,\sigma(W)\bigr)+E\bigl(\Var(\wt D(0,rZ|W|^\gamma)\,\big|\,\sigma(W))\bigr)
\\
\le
4\Var\bigl(E(\wt D(0,r^\gamma Z|W|^\gamma)\,\big|\,\sigma(W))\bigr)+2E\bigl(\Var(\wt D(0,r^\gamma Z|W|^\gamma)\,\big|\,\sigma(W))\bigr)
\\
+4E\bigl(\wt D(0,r^\gamma W)\bigr)^2-E\bigl(\wt D(0,rW)\bigr)^2+F_0\,.
\quad
\end{multline}
Abbreviating
\begin{equation}
\label{E:3.14ua}
\begin{aligned}
A_n&:=\frac1{4^n}\Var\Bigl(E\bigl(\wt D(0,r^{\gamma^{-n}}Z|W|^\gamma)\,\big|\,\sigma(W)\bigr)\Bigr)
\\
B_n&:=\frac1{4^n}E\Bigl(\Var\bigl(\wt D(0,r^{\gamma^{-n}}Z|W|^\gamma)\,\big|\,\sigma(W)\bigr)\Bigr)
\\
C_n&:=\frac1{4^n}E\bigl(\wt D(0,r^{\gamma^{-n}}W)\bigr)^2
\end{aligned}
\end{equation}
the inequality \eqref{E:3.13} gives
\begin{equation}
\label{E:3.14th}
A_n+B_n+C_n\le A_{n-1}+\frac12 B_{n-1}+C_{n-1}+\frac{F_n}{4^n}\,,
\end{equation}
where $F_n(r):=F_0(r^{\gamma^{-n}})$.
Iterating shows
\begin{equation}
\label{E:3.15th}
A_n+ \frac12 B_n+C_n\le A_0+\frac12B_0+C_0-\frac12\sum_{k=1}^n B_k + \sum_{k=1}^n \frac{F_k}{4^k}.
\end{equation}
Thanks to \eqref{E:2.6} and \eqref{E:3.4th} we have $\sup_{n\ge1}F_n/2^n<\infty$.
Since $A_n,B_n,C_n\ge0$ and~$F_n/4^n$ is summable on $n\ge0$, the sum of~$B_k$ must remain bounded uniformly in $n$.
\end{proofsect}

As a direct consequence we get:

\begin{corollary}
\label{cor-3.5}
Assume~$\wt D$ and~$W$ are independent with distributions as above. Then
\begin{equation}
\sup_{r\in[\texte^\gamma,\texte)}\,\sup_{n\ge1} E\biggl(\Bigl(\frac{\wt D(0,r^{\gamma^{-n}}W)}{2^n}\Bigr)^2\biggr)<\infty.
\end{equation}
\end{corollary}

\begin{proofsect}{Proof}
In the notation of the previous proof, the expectation equals~$A_n+B_n+C_n$ which is bounded uniformly in~$n$ thanks to \eqref{E:3.15th}. As $\sup_{n\ge1}(F_n(r)/2^n)$ is bounded uniformly on compact intervals of~$r$, the expectation is bounded also uniformly in~$r$ on the stated interval. 
\end{proofsect}

We are now ready to give:

\begin{proofsect}{Proof of Proposition~\ref{prop-2.5}}
Consider again the independent copies~$\wt D'$ and~$Z'$ of the quantities~$\wt D$ and~$Z$, respectively. 
Formula \eqref{E:1.14th} then reads
\begin{equation}
\sum_{n=1}^\infty E\Biggl[\biggl(\frac{\wt D'(0,r^{\gamma^{-n}}Z'|W|^\gamma)}{2^n}-\frac{\wt D(0,r^{\gamma^{-n}} Z|W|^\gamma)}{2^n}\biggr)^2\Biggr]<\infty\,.
\end{equation}
Pick a compact set $U\subset\mathbb{R}^d\setminus\{0\}$ with non-empty interior, denote its Lebesgue measure by~$|U|$ and let~$\epsilon\in(0,1)$. From the fact that $Z$ has a continuous nonvanishing density~$f$, there is a constant $c=c(U,\epsilon)>0$ such that
\begin{equation}
z|w|^\gamma\in U\quad\&\quad|w|<1/\epsilon\quad\Rightarrow\quad
f(z)|w|^{-d\gamma}\ge c\frac1{|U|}\,.
\end{equation}
Restricting the expectation to the event~$\{Z|W|^\gamma\in U\}\cup\{|W|<1/\epsilon\}$, this bound permits us to change variables from~$z$ to $x:=z|w|^\gamma$ and conclude that for~$X$ uniform on~$U$, and independent of all other random objects, we have
\begin{equation}
\sum_{n=1}^\infty E\Biggl[\biggl(\frac{\wt D'(0,r^{\gamma^{-n}}Z'|W|^\gamma)}{2^n}-\frac{\wt D(0,r^{\gamma^{-n}} X)}{2^n}\biggr)^2\,\Bigg|\,|W|<1/\epsilon\Biggr]<\infty
\end{equation}
where we also used that $P(|W|<1/\epsilon)>0$ for~$\epsilon\in(0,1)$. Using Jensen's inequality, we can now pass the expectation over~$\wt D$, $Z'$ and~$W$ inside the square to get
\begin{equation}
\label{correct}
\sum_{n=1}^\infty E\Biggl[\biggl(E\biggl[\frac{ \wt D(0,r^{\gamma^{-n}} Z|W|^\gamma)}{2^n}\,\bigg|\,|W|<1/\epsilon\biggr]-\frac{\wt D(0,r^{\gamma^{-n}} X)}{2^n}\biggr)^2\Biggr]<\infty
\end{equation}
By the Monotone Convergence Theorem, this implies
\begin{equation}
\label{convergence}
\frac{\wt D(0,r^{\gamma^{-n}} X)}{2^n}-E\biggl[\frac{ \wt D(0,r^{\gamma^{-n}} Z|W|^\gamma)}{2^n}\,\bigg|\,|W|<1/\epsilon\biggr]
\,\underset{n\to\infty}\longrightarrow\,0,\text{a.s.}
\end{equation}
with the exceptional set not depending on $\epsilon$.

Let~$\tilde c$ denote the quantity in Corollary~\ref{cor-3.5}. Denoting $q_\epsilon:=P(|W|\ge1/\epsilon)$, from Cauchy-Schwarz we have
\begin{equation}
\Biggl|
(1-q_\epsilon)E\biggl[\frac{ \wt D(0,r^{\gamma^{-n}} Z|W|^\gamma)}{2^n}\,\bigg|\,|W|<1/\epsilon\biggr]
-E\biggl[\frac{ \wt D(0,r^{\gamma^{-n}} W)}{2^n}\biggr]\Biggr|
\le \sqrt{\tilde cq_\epsilon} \label{eps}
\end{equation}
As $q_\epsilon\to0$ when $\epsilon\downarrow0$, we thus get \eqref{E:3.6th} for Lebesgue a.e. $x\in U$. Since~$U$ was arbitrary (compact), the same applies to a.e.~$x\in\R^d$. 
\end{proofsect}

\begin{remark}
The reader may wonder why the passage through an a.s.\ limit for Lebesgue a.e.~$x$ has been used instead of trying to prove the a.s.\ convergence of $X_n:=2^{-n}\wt D(0,r^{\gamma^{-n}}W)$ directly. (The convergence $X_n\to L(r)$ a.s.\ does hold by \eqref{E:3.6th} and the fact that~$W$ has a density w.r.t.\ the Lebesgue measure.) This is because Lemma~\ref{lemma-2.6} only controls the conditional variances of~$X_n$ given~$W$, and not the full variances $\Var(X_n)$. We will in fact show $\Var(X_n)\to0$ in the proof of Proposition~\ref{prop-3.6}, but that only with the help of~\eqref{E:3.6th}. 
\end{remark}

\subsection{Full limit for the restricted distance}
We now proceed to extend the limit from multiples of the argument by terms from $\{r^{\gamma^{-n}}\colon n\ge0\}$ to multiples ranging continuously through positive reals. However, for reasons described after Theorem~\ref{thm-2}, such a limit can generally be claimed only in probability. It will also suffice to show this for~$x$ replaced by the random variable~$W$. This is the content of:

\begin{proposition}
\label{prop-3.6}
Suppose~$\wt D$ and~$W$ are independent with distributions as above. Then
\begin{equation}
\label{E:3.24ua}
\frac{\wt D(0,r W)}{L(r)}\,\underset{r\to\infty}\longrightarrow\,1\quad\text{in probability}.
\end{equation}
\end{proposition}

As we will see, a key point in proving Proposition~\ref{prop-3.6} is:

\begin{lemma}
\label{lemma-3.7}
The identity $L(r)=2L(r^\gamma)$ holds for all~$r\ge1$ and $t\mapsto L(\texte^t)$ is convex on~$[0,\infty)$. In particular,~$r\mapsto L(r)$ is continuous, strictly increasing on~$[1,\infty)$. The function
\begin{equation}
\phi(r):=L(r)(\log r)^{-\Delta},\qquad r>1,
\end{equation}
obeys the conditions stated in Theorem~\ref{thm-2}.
\end{lemma}

\begin{proofsect}{Proof}
First, $L(r)=2L(r^\gamma)$ is a consequence of the limit definition of~$L$ in Lemma~\ref{lemma-2.6}.
Let~$\gamma_1,\gamma_2$ be such that $0<\gamma_1,\gamma_2<\frac{1+\gamma}2$ and~$d\gamma_1+d\gamma_2=s$. Plugging~$r^{\gamma^{-n}}x$ for~$x$ in \eqref{E:2.7} yields
\begin{equation}
\wt D(0,r^{\gamma^{-n}}x)\,\overset{\text{\rm law}}\le\, \wt D\bigl(0,r^{\gamma_1\gamma^{-n}}|x|^\gamma Z\bigr)+\wt D'\bigl(0,r^{\gamma_2\gamma^{-n}}|x|^{\gamma_2}Z'\bigr)+1+r^{\gamma^{-n}}|x|1_{A(r^{\gamma^{-n}}x)}.
\end{equation}
Dividing the expression by~$2^n$, applying  
Proposition~\ref{prop-2.5} and noting that, by \eqref{E:2.6}, the last two terms tend to zero in probability as~$n\to\infty$ gives
\begin{equation}
L(r)\le L(r^{\gamma_1})+L(r^{\gamma_2}).
\end{equation}
Now set $\texte^{t_1}:=r^{\gamma_1/\gamma}$ and $\texte^{t_2}:=r^{\gamma_2/\gamma}$ and observe that then $r=\texte^{\frac12(t_1+t_2)}$. The identity $2L(r^\gamma)=L(r)$ and the fact that the constraints on~$\gamma_1,\gamma_2$ will be satisfied if~$|\gamma_1-\gamma_2|<1-\gamma$ then imply
\begin{equation}
\forall t_1,t_2\ge0\colon\quad
0<\frac{|t_1-t_2|}{t_1+t_2}<\frac{1-\gamma}{2\gamma}
\quad\Rightarrow\quad
L\bigl(\texte^{\frac12(t_1+t_2)}\bigr)\le\frac{L(\texte^{t_1})+L(\texte^{t_2})}2\,,
\end{equation}
i.e., a local mid-point convexity of~$t\mapsto L(\texte^t)$. The upper semicontinuity of~$L$ from Lemma~\ref{lemma-2.6} then implies continuity of~$r\mapsto L(r)$ on~$[1,\infty)$, and subsequently also the convexity of $t\mapsto L(\texte^t)$ on~$[0,\infty)$. The strict monotonicity arises from convexity and the fact that $L(1)=0$ while~$L(r)>0$ for~$r>1$, by Lemma~\ref{lemma-3.2}. The conditions for~$\phi$ in~Theorem~\ref{thm-2}  are checked directly.
\end{proofsect}

We will also need a uniform bound on third moments of~$2^{-n}\wt D(0,r^{\gamma^{-n}}W)$:

\begin{lemma}
\label{lemma-3.9}
Assume~$\wt D$ and~$W$ are independent with distributions as above. Then
\begin{equation}
\sup_{r\in[\texte^\gamma,\texte)}\,\sup_{n\ge1}\,\,E\biggl(\Bigl(\frac{\wt D(0,r^{\gamma^{-n}}W)}{2^n}\Bigr)^3\biggr)<\infty.
\end{equation}
\end{lemma}

\begin{proofsect}{Proof}
Consider the setup of Proposition~\ref{prop-2.2} with~$\gamma_1=\gamma_2:=\gamma$. Taking the third power of both sides of \eqref{E:2.7} and setting with~$x:=rW$ yields
\begin{equation}
E\bigl(\wt D(0,rW)^3\bigr)\le E\Bigl(\bigl[\wt D(0,r^\gamma|W|^\gamma Z)+\wt D'(0,r^\gamma|W|^\gamma Z)\bigr]^3\Bigr)+G_0(r)\,,
\end{equation}
where
\begin{multline}
G_0(r):=3E\Bigl(\bigl[\wt D(0,r^\gamma|W|^\gamma Z)+\wt D'(0,r^\gamma|W|^\gamma Z)\bigr]^2\Bigr)
\\
+3E\Bigl(\wt D(0,r^\gamma|W|^\gamma Z)+\wt D'(0,r^\gamma|W|^\gamma Z)\Bigr)+1+3 E\Bigl(\bigl[1+2r^\gamma|W|^\gamma\bigr]^2r|W|1_{A(rW)}\Bigr)
\\+3E\Bigl(\bigl[1+2r^\gamma|W|^\gamma\bigr]r^2|W|^21_{A(rW)}\Bigr)+E\Bigl(r^3|W|^31_{A(rW)}\Bigr)
\end{multline}
With the help of H\"older's inequality and the fact that $Z|W|^\gamma\laweq W$, we then get
\begin{equation}
E\bigl(\wt D(0,r^{\gamma^{-n}}W)^3\bigr)\le 8 E\bigl(\wt D(0,r^{\gamma^{-n+1}}W)^3\bigr)+G_n(r)\,,
\end{equation}
where, as before, $G_n(r):=G_0(r^{\gamma^{-n}})$. Corollary~\ref{cor-2.5} and \eqref{E:2.6} ensure that $G_n(r)/4^n$ is bounded uniformly in~$n\ge1$ and~$r\in[\texte^\gamma,\texte)$. The claim follows.
\end{proofsect}

We are now ready to give:

\begin{proofsect}{Proof of Proposition~\ref{prop-3.6}}
Abbreviate $X_n:=2^{-n}\wt D(0,r^{\gamma^{-n}}W)$. By Lemma~\ref{lemma-3.2}, $E(X_n)\to L(r)$. Lemma~\ref{lemma-3.9} and the almost sure convergence in Proposition~\ref{prop-2.5} in turn show $E(X_n^2)\to L(r)^2$. Using the quantities from \eqref{E:3.14ua}, we can alternatively write $E(X_n^2)=A_n+B_n+C_n$. Since~$C_n=[EX_n]^2$, the above shows $C_n\to L(r)^2$ and so~$A_n+B_n\to0$ for each~$r\ge1$. 

We claim that the convergence of the moments is uniform in~$r$ on compact subsets of~$[1,\infty)$. Starting with the former, observe that \eqref{E:3.14th} in fact shows that
\begin{equation}
n\mapsto A_n+B_n+C_n-\sum_{k=n+1}^\infty\frac{F_k}{4^k}
\end{equation}
is non-increasing, with the sum absolutely convergent. The functions $A_n,B_n,C_n,F_n$ are continuous and so is the limit~$L(r)^2$, by Lemma~\ref{lemma-3.7}. Dini's Theorem then ensures that the convergence $A_n+B_n\to0$ is indeed uniform on compact sets of~$r$. 

The argument for the first moments is similar; the proof of Lemma~\ref{lemma-3.2} shows that, for some constant~$c>0$, the sequence $E(X_n)+c2^{-n}$ decreases to~$L(r)$ pointwise. Since both the sequence and the limit are continuous, Dini's Theorem again implies local uniformity.

As $A_n+B_n = \Var(X_n)$, we have $\Var(X_n)\to0$ locally uniformly in~$r$. In light of the similar uniformity of~$E(X_n)\to L(r)$, for each~$\epsilon>0$ there is~$n_0\ge1$ such that
\begin{equation}
\sup_{r\in[\texte^\gamma,\texte)}\,\sup_{n\ge n_0}\,\,E\biggl(\Bigl|\frac{\wt D(0,r^{\gamma^n}W)}{2^n}-L(r)\Bigr|^2\biggr)<\epsilon.
\end{equation}
Since~$L(r^{\gamma^{-n}})=2^n L(r)$, dividing the expression by $L(r)$ we get
\begin{equation}
\sup_{r\ge\texte^{\gamma^{1-n_0}}}\,\,E\biggl(\Bigl|\frac{\wt D(0,rW)}{L(r)}-1\Bigr|^2\biggr)<\frac{\epsilon}{L(\texte^\gamma)},
\end{equation}
where we also used that~$L(r)\ge L(\texte^\gamma)$ for all $r\in[\texte^\gamma,\texte)$. This implies $\wt D(0,rW)/L(r)\to1$ in~$L^2$ and thus in probability.
\end{proofsect}

\begin{remark}
\label{rem-3.10}
With some extra work, we could show that the limit in \eqref{E:3.24ua} also exists in probability conditional on~$W$. As~$W$ is continuously distributed with support~$\R^d$, we could then replace~$W$ by Lebesgue a.e.~$x\in\R^d$ and, finally, use monotonicity arguments to extend to conclusion to all non-zero~$x\in\R^d$. However, the same arguments will (have to) be applied to the full distance treated in the next section and so we refrain from making them here.
\end{remark}

\section{Full distance scaling}
\label{sec4}\nopagebreak\noindent
We are now ready to return to the full distance~$D(x,y)$ associated with long-range percolation on~$\R^d$ and prove its asymptotic stated in Theorem~\ref{thm-2}. We begin by extending the conclusions of Proposition~\ref{prop-3.6} to the full distance:

\begin{proposition}
\label{prop-4.1}
Suppose~$D$ and~$W$ are independent with distributions as above. Then
\begin{equation}
\label{E:3.24uab}
\frac{D(0,r W)}{L(r)}\,\underset{r\to\infty}\longrightarrow\,1\quad\text{in probability}.
\end{equation}
\end{proposition}

For this, we will need to expand the notion of the restricted distance to a whole family of ``distances'' $D_k$ indexed by $k\in\{0,1,\dots\}$ as follows. Abbreviating~$\tilde\gamma:=\frac12(1+\gamma)$, we set
\begin{equation}
\label{E:2.1a}
\wt D_k(x,y):=\inf\Biggl\{n+\sum_{i=0}^n|x_{i+1}-y_i|\,\colon\,
\begin{aligned}
&n\ge0,\,\{(x_i,y_i)\colon i=1,\dots,n\}\subset\II_\beta
\\*[-1mm]
&x_i,y_i\in B\bigl(x,2|x-y|^{\tilde\gamma^{-k}}\bigr)\,\forall i=1,\dots,n
\end{aligned}
\Biggr\}\,,
\end{equation}
where, as before, we set~$y_0:=x$ and~$x_{n+1}:=y$. We have
\begin{equation}
D(x,y)\le\dots\le \wt D_{k+1}(x,y)\le \wt D_k(x,y)\le\dots\le \wt D_1(x,y)\le\wt D_0(x,y)=\wt D(x,y)\,.
\end{equation}
Our first observation is:

\begin{lemma}
\label{lemma-4.2}
Let~$W$ be independent of the distances~$\wt D_k$ and~$D$. There is $k\in\N$ such that
\begin{equation}
\label{E:4.4}
\lim_{r\to\infty}\,P\bigl(\,\wt D_k(0,rW)\ne D(0,rW)\bigr)=0.
\end{equation}
\end{lemma}

\begin{proofsect}{Proof}
Pick~$x\in\R^d$ and let~$c$ denote the diameter of~$[0,1)^d$ in~$|\cdot|$-norm. Note that, as soon as~$r$ is sufficiently large, on $\{\wt D_k(0,rx)\ne D(0,rx)\}$ there must be a point
\begin{equation}
y\in \Z^d\smallsetminus\Bigl[ B\bigl(0,(r|x|)^{\tilde\gamma^{-k}}\bigr)\cup B\bigl(x,(r|x|)^{\tilde\gamma^{-k}}\bigr)\Bigr]
\end{equation}
for which $D(0,y)\le \wt D(0,rx)+c$ and $D(rx,y)\le \wt D(0,rx)+c$ occur using disjoint collections of edges in the underlying sample of the Poisson process. Given any~$C>0$ and assuming that~$r$ is so large that~$C(\log r)^\Delta>2c$, the van den Berg-Kesten inequality then shows
\begin{multline}
\quad
P\Bigl(\,\wt D_k(0,rx)\ne D(0,rx),\,\wt D(0,rx)\le C(\log r)^\Delta\Bigr)
\\
\le\sum_{\begin{subarray}{c}
y\in\Z^d\\|y|\wedge|\lfloor rx\rfloor-y|\ge(r|x|)^{\tilde\gamma^{-k}}
\end{subarray}}
P\bigl(D(0,y)\le 2C(\log r)^\Delta\bigr)P\bigl(D(\lfloor rx\rfloor,y)\le 2C(\log r)^\Delta\bigr)\,.
\quad
\end{multline}
where~$a\wedge b:=\min\{a,b\}$ and where $\lfloor rx\rfloor$ is the closest point on~$\Z^d$ to~$rx$.
Our aim is to show that the sum vanishes as~$r\to\infty$ once~$k$ is large enough.

In light of the domination bound in Lemma~\ref{lemma-3th}, we can replace the continuum distance~$D$ by the discrete distance~$\Ddis$ at the cost of changing~$C$ and~$\beta$ by multiplicative constants. We may thus estimate the above sum for the model on~$\Z^d$ instead, writing temporarily just~$x$ for~$\lfloor rx\rfloor$ and~$n$ for~$2C(\log r)^\Delta$. The bound in Theorem~\ref{thm-4th} then shows
\begin{multline}
\qquad
\sum_{\begin{subarray}{c}
y\in\Z^d\\|y|\wedge|y-x|\ge|x|^{\tilde\gamma^{-k}}
\end{subarray}}\,
P\bigl(\Ddis(0,y)\le n\bigr)P\bigl(\Ddis(y,x)\le n\bigr)
\\*[-5mm]
\le c_1^2 \texte^{2c_2n^{1/\Delta}}
\sum_{\begin{subarray}{c}
y\in\Z^d\\|y|\wedge|y-x|\ge|x|^{\tilde\gamma^{-k}}
\end{subarray}}
\frac1{|y|^s|x-y|^s}
\le \tilde c_1\frac{\texte^{2 c_2'n^{1/\Delta}}}{|x|^{s\tilde\gamma^{-k}}}
\qquad
\end{multline}
for some~$\tilde c_1\in(0,\infty)$ independent of~$x$.
Returning to the continuum problem with $n:=C(\log r)^{\Delta}$, there is thus a constant~$\tilde c\in(0,\infty)$ such that, for some~$c_2'$ proportional to~$c_2$,
\begin{equation}
P\Bigl(\,\wt D_k(0,rx)\ne D(0,rx),\,\wt D(0,rx)\le C(\log r)^\Delta\Bigr)
\le \tilde c\frac{|x|^{-s\tilde\gamma^{-k}}}{r^{s\tilde\gamma^{-k}-2 c_2'C^{1/\Delta}}}\,.
\end{equation}
The exponent of~$r$ in the denominator is positive once~$k$ is taken sufficiently large (depending only on~$C$). Plugging in~$x:=W$, choosing~$C>\max\phi$ for~$\phi$ as in Lemma~\ref{lemma-3.7}, adjusting~$k$ accordingly and invoking Proposition~\ref{prop-3.6}, we get \eqref{E:4.4} on the event~$\{|W|>\epsilon\}$, for any~$\epsilon>0$. But~$W$ is continuously distributed and so the claim follows by noting that $P(|W|\le\epsilon)\to0$ as~$\epsilon\downarrow0$. 
\end{proofsect}

Next we observe:

\begin{lemma}
\label{lemma-4.4}
Let~$k\in\N$ and suppose~$W$ and~$\wt D_k$ are independent with distributions as above. Then for every~$k\ge0$,
\begin{equation}
\label{E:4.9}
\liminf_{r\to\infty}\,\,\frac{E \wt D_k(0,rW)}{L(r)}\ge1.
\end{equation}
\end{lemma}

The proof will be based on perturbations of the underlying model in~$\beta$. For this reason, let~$L_\beta(r)$ henceforth mark the explicit dependence of the limit in Lemma~\ref{lemma-3.2} on~$\beta$ and define, as before, $\phi_\beta(r):=L_\beta(r)/(\log r)^\Delta$. We note one useful fact:

\begin{lemma}
\label{lemma-4.5}
The function~$(\beta,r)\mapsto\phi_\beta(r)$ is jointly continuous on~$(0,\infty)\times(1,\infty)$.
\end{lemma}

\begin{proofsect}{Proof}
From \eqref{E:1.7w} and the existence of the limit in Lemma~\ref{lemma-3.2} we get
\begin{equation}
\forall a\ge1\colon\quad L_\beta(r)\le L_{a^{s-2d}\beta}(r)\le a L_\beta(r).
\end{equation}
The continuity of~$\beta\mapsto L_\beta(r)$ for each~$r\ge1$ is then readily inferred. The continuity and monotonicity of $r\mapsto L_\beta(r)$ then yields the joint continuity of $(\beta,r)\mapsto L_\beta(r)$ on~$(0,\infty)\times[1,\infty)$. The claim follows by applying the definition of~$\phi_\beta$.
\end{proofsect}

\begin{proofsect}{Proof of Lemma~\ref{lemma-4.4}}
The argument hinges on a subadditivity bound of the kind derived in Proposition~\ref{prop-2.2} which links expectations of~$\wt D_k$ and~$\wt D_{k+1}$ albeit at slightly different values of~$\beta$. The proof of this bound follows closely that of the above proposition, although it is simpler as here we can efficiently use additivity of Poisson processes.  

Fix~$\beta>0$ and let~$W$ be the random variable associated with parameters~$\beta$ and~$\eta:=1$ as defined in Lemma~\ref{lemma-3.1}. Let~$\epsilon\in(0,1/2)$. Writing~$E_\beta$ for the expectation with respect to the point process~$\II_\beta$ with intensity measure~$\mu_{s,\beta}$, we will show later that, for some~$c=c(\beta,\epsilon)\in(0,\infty)$,
\begin{equation}
\label{E:4.11}
E_\beta\bigl[ D_k(0,r\epsilon^{-\frac1{2d}}W)\bigr]\le 2 E_{\beta(1-2\epsilon)}\bigl[D_{k+1}(0,r^\gamma \epsilon^{-\frac1{2d}}W)\bigr]+c
\end{equation}
holds for all~$r\ge1$. This is sufficient to prove the claim by induction. Indeed, the factors~$\epsilon^{-\frac1{2d}}$ can seamlessly be absorbed into~$r$ by noting that $L_\beta(ar)/L(r)\to1$ as~$r\to\infty$ for any~$a>0$ thanks to the continuity of~$r\mapsto\phi_\beta(r)$. Assuming \eqref{E:4.9} for some~$k\in\N$, then dividing \eqref{E:4.11} by~$L_\beta(r)=2L_\beta(r^\gamma)$ and relabeling~$\beta(1-2\epsilon)$ for~$\beta$ yields
\begin{equation}
\liminf_{r\to\infty} \frac{E_\beta\bigl[ D_{k+1}(0,rW)\bigr]}{L_\beta(r)}
\ge \inf_{r\in[\texte^{\gamma},\texte)}\,\frac{\phi_{\beta/(1-2\epsilon)}(r)}{\phi_\beta(r)}
\end{equation}
Taking $\epsilon\downarrow0$ and invoking the continuity from Lemma~\ref{lemma-4.5}, we then get \eqref{E:4.9} for~$k+1$ as well. Since Lemma~\ref{lemma-3.2} ensures \eqref{E:4.9} for~$k:=0$, we get it for all~$k\ge0$.
 
It remains to prove \eqref{E:4.11}. Abbreviate
\begin{equation}
\beta':=2\epsilon\beta\quad\text{and}\quad\beta'':=(1-2\epsilon)\beta.
\end{equation}
A sample~$\II_\beta$ of the Poisson process with intensity measure~$\mu_{s,\beta}$ can then be written as the union $\II_{\beta'}\cup\II_{\beta''}$ of two independent processes with intensities $\mu_{s,\beta'}$ and $\mu_{s,\beta''}$, respectively. Fix~$x\in\R^d$. Following the proof of Proposition~\ref{prop-2.2} (with~$\eta$ there set to $1/2$), under the condition
\begin{equation}
\label{E:4.12}
1+2|x|^{\tilde\gamma-1}\le 2^{1/s}
\end{equation}
we can further decompose~$\II_{\beta'}$ into the union of independent processes~$\II'_{\beta'}$ and~$\II''_{\beta'}$, with their respective intensity measures given by
\begin{equation}
\mu'_{s,\beta'}(\textd \tilde x\,\textd \tilde y):=\epsilon\beta\1_{\{|\tilde x|_2<|\tilde y|_2\}}\1_{\bigl\{|\tilde x|\vee|x-\tilde y|<|x|^{\tilde\gamma}\bigr\}}\frac{\textd\tilde x\,\textd\tilde y}{|x|^s}
\end{equation}
and $\mu''_{s,\beta'}:=\mu_{s,\beta'}-\mu'_{s,\beta'}$.
(The condition \eqref{E:4.12} ensures that~$\mu''_{s,\beta'}$ is a positive measure.) We also introduce an auxiliary independent process~$\II_{\beta'}'''$ with intensity measure
\begin{equation}
\mu'''_{s,\beta'}(\textd \tilde x\,\textd \tilde y):=\epsilon\beta\Bigl(1-\1_{\{|\tilde x|_2<|\tilde y|_2\}}\1_{\bigl\{|\tilde x|\vee|x-\tilde y|<|x|^{\tilde\gamma}\bigr\}}\Bigr)\frac{\textd\tilde x\,\textd\tilde y}{|x|^s}\,.
\end{equation}
As is directly checked, $\II'_{\beta'}\cup\II'''_{\beta'}$ is a homogenous Poisson process with intensity~$\epsilon\beta|x|^{-s}$.

Now define a pair of random variables~$(X,Y)$ as the minimizer of $f_x(\tilde x,\tilde y):=|\tilde x|^{2d}+|x-\tilde y|^{2d}$ among all points of~$\II'_{\beta'}\cup\II'''_{\beta'}$. Set
\begin{equation}
Z:=\epsilon^{\frac1{2d}}|x|^{-\gamma}X\quad\text{and}\quad Z':=\epsilon^{\frac1{2d}}|x|^{-\gamma}(x-Y)
\end{equation}
and note that, by the calculation in \eqref{E:2.13} and a simple scaling arguement, $Z,Z'$ are i.i.d.\ with common law \eqref{E:2.6a} for~$\eta:=1$.  Given~$k\in\N$, let~$D_k(0,x)$ be defined using the full process~$\II_\beta$ and let~$D''_{k+1}(\cdot,\cdot)$ be defined using the process~$\II_{\beta''}$. We now claim
\begin{equation}
\label{E:4.17}
D_k(0,x)\le D''_{k+1}\bigl(0,\epsilon^{-\frac1{2d}}|x|^\gamma Z\bigr)+D''_{k+1}\bigl(x,x+\epsilon^{-\frac1{2d}}|x|^\gamma Z'\bigr)+1+|x|1_{A'(x)}
\end{equation}
where we set
\begin{equation}
A'(x):=\bigl\{|Z|\vee|Z'|\ge\tfrac12\epsilon^{\frac1{2d}}|x|^{\gamma-\tilde\gamma}\bigr\}.
\end{equation}
whenever~$|x|$ is so large that \eqref{E:4.12} and
\begin{equation}
\label{E:4.19}
|x|^{1-1/\tilde\gamma^k}+2^{1-1/\tilde\gamma^{k+1}}\le 2,\qquad k\ge0,
\end{equation}
hold, and put $A'(x)$ to the whole sample space otherwise.
To see this we note that, on $A'(x)$ the inequality follows from $D_k(0,x)\le|x|$ and so we just need to prove this on~$A'(x)^\cc$. Here we observe that $|X|\vee|x-Y|\le\frac12|x|^{\tilde\gamma}$ and so~$(X,Y)\in\II_{\beta'}$. A path minimizing $D''_{k+1}(0,X)$ will then lie in $B(0,2|X|^{\tilde\gamma^{-k-1}})\subseteq B(0,2|x|^{\tilde\gamma^{-k}})$ while the path minimizing $D''_{k+1}(x,Y)$ will lie in
\begin{equation}
B\bigl(x,2|x-Y|^{\tilde\gamma^{-k-1}}\bigr)\subseteq B\bigl(x,2^{1-\tilde\gamma^{-k-1}}|x|^{\tilde\gamma^{-k}}\bigr)\subseteq B\bigl(0,2|x|^{\tilde\gamma^{-k}}\bigr)\,,
\end{equation}
where the last inclusion is inferred from \eqref{E:4.19}. The concatenation of these paths with edge~$(X,Y)$ then produces a path entering the infimum defining~$D_k(0,x)$. Hence \eqref{E:4.17} follows.

Noting that the probability of~$A'(x)$ decays stretched-exponentially with~$|x|$, plugging~$W$ on both sides of \eqref{E:4.17}, taking expectation and using that~$|W|^\gamma Z\laweq W$ then yields \eqref{E:4.11}.
\end{proofsect}

Armed with the above lemmas, we can now give:

\begin{proofsect}{Proof of Proposition~\ref{prop-4.1}}
Abbreviate~$X(r):=\wt D_k(0,rW)/L(r)$. Corollary~\ref{cor-3.5} and $D(0,rW)\le \wt D(0,rW)$ show $\sup_{r\ge\texte^\gamma}E[X(r)^2]<\infty$ and so~$\{X(r)\colon r\ge\texte^\gamma\}$ is uniformly integrable. Proposition~\ref{prop-3.6} in turn implies~$P(X(r)>1+\epsilon)\to0$ as~$r\to\infty$ for every~$\epsilon>0$ and so we have $E[X(r)1_{\{X(r)>1+\epsilon\}}]\to0$ as well. Lemma~\ref{lemma-4.4} then gives~$E[X(r)]\to1$. Since~$X(r)\ge0$, it follows that the mass of~$X(r)$ must asymptotically concentrate at~$1$. This proves the claim for~$\wt D_k(0,rW)$; Lemma~\ref{lemma-4.2} then extends it to~$D(0,rW)$.
\end{proofsect}

This makes us finally ready to complete the proof of our main results:

\begin{proofsect}{Proof of Theorem~\ref{thm-2}}
The definition and properties of function~$\phi$ have already been established, so we just have to prove the limit claim \eqref{E:1.8}. This will be derived from Proposition~\ref{prop-4.1} and some perturbation arguments. Write~$P_\beta$ for the law of the Poisson process with intensity~$\mu_{s,\beta}$. Fix~$x\in\R^d\smallsetminus\{0\}$ and note that, by the stochastic domination bounds in Lemmas~\ref{lemma-4.3}--\ref{lemma-4.3a}, for each~$\epsilon$ there is~$\delta$ such that for all~$y\in\R^d$ and all~$t>0$,
\begin{equation}
|y-x|<\delta|x|\quad\Rightarrow\quad P_\beta\bigl(D(0,x)\le t\bigr)\le P_{\beta(1+\epsilon)}\bigl(D(0,y)\le(1+\epsilon)t\bigr)\,.
\end{equation}
Let~$W$ be independent of~$D$ with the distribution as above and pick any~$\zeta\in(0,1)$. Noting that~$P(|W-x|<\delta|x|)>0$, for the above~$\epsilon$ and~$\delta$ we then get
\begin{equation}
\label{E:4.23}
\begin{aligned}
P_\beta\bigl(D(0,rx)\le&\,(1-\zeta)L_\beta(r)\bigr)
\\
&=P_\beta\Bigl(D(0,rx)\le(1-\zeta)L_\beta(r)\,\Big|\,|W-x|<\delta|x|\Bigr)
\\
&\le P_{\beta(1+\epsilon)}\Bigl(D(0,rW)\le(1-\zeta)(1+\epsilon)L_\beta(r)\,\Big|\,|W-x|<\delta|x|\Bigr)
\end{aligned}
\end{equation}
Lemma~\ref{lemma-4.5} permits us to pick~$\epsilon$ so small that
\begin{equation}
(1-\zeta)(1+\epsilon)\inf_{r\in[\texte^\gamma,\texte)}\frac{\phi_\beta(r)}{\phi_{\beta(1+\epsilon)}(r)}<1-\epsilon.
\end{equation}
The right-hand side of \eqref{E:4.23} then tends to zero by Proposition~\ref{prop-4.1}. The argument for the other bound is completely analogous and so we omit it.
\end{proofsect}

\begin{proofsect}{Proof of Theorem~\ref{thm-1}}
The lower bound was already shown in Corollary~\ref{cor-2.5}. For the the upper bound we first use Theorem~\ref{thm-2} and the comparisons in Lemma~\ref{lemma-3th} to prove the claim for~$x$ of the form~$x:=r\texte_i$, where~$\texte_i$ is one of the coordinate vectors. Then we use the triangle inequality for~$\Ddis$ to get the full limit as~$|x|\to\infty$.
\end{proofsect}

\section*{Acknowledgments}
\nopagebreak\noindent
This research has been partially supported by NSF grant DMS-1407558 and GA\v CR project P201/16-15238S. 
Part of this work has appeared as part of the second author's PhD~thesis~\cite{Lin-thesis}.

\end{document}